\documentclass{article}    
\usepackage{amssymb, latexsym, amsmath, amsthm}

\newcommand{\bi} {\begin{itemize}}
\newcommand{\ei} {\end{itemize}}
\newcommand{\bea} {\begin{eqnarray}}
\newcommand{\eea} {\end{eqnarray}}
\newcommand{\be} {\begin{equation}}
\newcommand{\ee} {\end{equation}}
\newcommand{\bean} {\begin{eqnarray*}}
\newcommand{\eean} {\end{eqnarray*}}

\RequirePackage{epsfig}
\RequirePackage{amssymb}
\RequirePackage{natbib}
\RequirePackage{graphics}
\bibpunct{(}{)}{,}{a}{,}{,}

\begin{document}
\title{Bounded Point Evaluations For Rationally Multicyclic Subnormal Operators}
\author{Liming Yang}
\author{
Liming Yang \\
 Department of Mathematics \\
Virginia Polytechnic Institute and State  University \\
Blacksburg, VA 24061 \\
yliming@vt.edu
     }         
\date{}
\maketitle

\newtheorem{Theorem}{Theorem}
\newtheorem*{MTheorem}{Main Theorem}
\newtheorem*{TTheorem}{Thomson's Theorem}
\newtheorem*{Conjecture}{Conjecture}
\newtheorem{Corollary}{Corollary}
\newtheorem{Definition}{Definition}
\newtheorem{Assumption}{Assumption}
\newtheorem{Lemma}{Lemma}
\newtheorem{Problem}{Problem}
\newtheorem*{Example}{Example}
\newtheorem*{Remark}{Remark}
\newtheorem{KnownResult}{Known Result}
\newtheorem{Algorithm}{Algorithm}
\newtheorem{Property}{Property}
\newtheorem{Proposition}{Proposition}


\abstract{
Let $S$ be a pure bounded rationally multicyclic subnormal operator on a separable complex Hilbert space $\mathcal H$ and let $M_z$ be the minimal normal extension on a separable complex Hilbert space $\mathcal K$ containing $\mathcal H.$ Let $bpe(S)$ be the set of bounded point evaluations and let $abpe(S)$ be the set of analytic bounded point evaluations.
We show $abpe(S) = bpe(S) \cap Int(\sigma (S)).$ The result affirmatively answers a question asked by J. B. Conway concerning the equality of the interior of $bpe(S)$ and $abpe(S)$ for a rationally multicyclic subnormal operator $S.$ As a result, if $\lambda_0\in Int(\sigma (S))$ and $dim(ker(S-\lambda_0)^*) = N,$ where $N$ is the  minimal number of cyclic vectors for $S,$ then the range of $S-\lambda_0$ is closed, hence, $\lambda_0 \in \sigma (S) \setminus \sigma_e (S).$   
}

\section{Introduction}

Let $\mathcal H$ be a separable complex Hilbert space and let $\mathcal L(\mathcal H )$ be the space of bounded linear operators on $\mathcal H.$ An operator $S\in \mathcal L (\mathcal H)$ is subnormal if there exist a separable complex Hilbert space $\mathcal K$ containing $\mathcal H$ and a normal operator $M_z \in \mathcal L(\mathcal K)$ such that $M_z\mathcal H \subset \mathcal H$ and $S = M_z|_{\mathcal H}.$ By the spectral theorem of normal operators, we assume that 
\[
\ \mathcal K = \oplus _{i = 1}^m L^2 (\mu _i) \tag{1-1}
\]
where $\mu_1 >> \mu_2 >> ... >>\mu_m$ ($m$ may be $\infty$) are compactly supported finite positive measures on the complex plane
$\mathbb {C},$ and $M_z$ is multiplication by $z$ on $\mathcal K.$ 
For $H= (h_1,...,h_m)\in \mathcal K$ and $G= (g_1,...,g_m)\in \mathcal K,$ we define
 \[
 \ \left\langle H(z), G(z) \right\rangle = \sum _{i = 1}^m h_i(z) \overline{g_i(z)} \dfrac{d\mu_i}{d\mu _1}, ~ | H(z) |^2 = \left\langle H(z), H(z) \right\rangle . \tag{1-2}
 \]
The inner product of $H$ and $G$ in $\mathcal K$ is defined by
 \[
 \ (H,G) = \int \left\langle H(z), G(z) \right\rangle d\mu_1(z).\tag{1-3}
 \]
$M_z$ is the minimal normal extension if 
\[
 \ \mathcal K = clos\left( span (M_z^{*k}x: ~x\in \mathcal H,~ k\ge 0)\right ). \tag{1-4}
 \]
We will always assume that $M_z$ is the minimal normal extension of $S$ and $\mathcal K$ satisfies (1-1) and (1-4). For details about the functional model above and basic knowledge of subnormal operators, the reader shall consult Chapter II of the book \cite{conway}.

For $T \in \mathcal L(\mathcal H )$, we denote by $\sigma (T)$ the
spectrum of $T,$ $\sigma _e (T)$ the essential spectrum
of $T,$ $T^*$ its adjoint, $ker(T)$ its kernel, and $Ran(T)$ its range. For a
subset $A \subset \mathbb C,$ we set $Int(A)$ for its interior, $\bar A$ or $clos (A)$  for its closure, $A^c$ for its complement, and $\chi _A$ for its characteristic function. Let $\delta_{ij} =1$ when $i=j$ and $\delta_{ij} = 0$ when $i\ne j.$ For $\lambda\in \mathbb C$ and $\delta > 0,$ we set $B(\lambda, \delta) = \{z: |z - \lambda  | <\delta \}$ and $\mathbb D = B(0,1).$ Let $\mathcal{P}$ denote the set of polynomials in the complex variable $z.$ For a compact subset $K\subset \mathbb C,$ let $Rat(K)$ be the set of all rational functions with poles off $K.$

A subnormal operator $S$ on $\mathcal H$ is pure if for every non-zero invariant subspace $I$ of $S$ ($SI \subset I$), the operator $S|_I$ is not normal. 
For $F_1, F_2,..., F_N\in \mathcal H,$ let
 \[
 \ R^2(S | F_1,F_2,...,F_N) = clos\{r_1(S)F_1 + r_2(S)F_2 + ... + r_N(S)F_N\} \tag{1-5}
 \]
 in $\mathcal H,$ where $r_1,r_2,...,r_N \in Rat(\sigma(S))$ and let
 \[
 \ P^2(S | F_1,F_2,...,F_N) = clos\{p_1(S)F_1 + p_2(S)F_2 + ... + p_N(S)F_N\} \tag{1-6}
 \] 
in $\mathcal H,$ where $p_1,p_2,...,p_N \in \mathcal P.$
A subnormal operator $S$ on $\mathcal H$ is rantionally multicyclic ($N-$cyclic) if there are $N$ vectors $F_1,F_2,...,F_N\in \mathcal H$ such that 
 \[
 \ \mathcal H = R^2(S | F_1,F_2,...,F_N)
 \]
and for any $G_1,...,G_{N-1}\in \mathcal H,$
\[
 \ \mathcal H \ne R^2(S | G_1,G_2,...,G_{N-1}).
 \]
$S$ is multicyclic ($N-$cyclic) if 
 \[
 \ \mathcal H = P^2(S | F_1,F_2,...,F_N)
 \]
and for any $G_1,...,G_{N-1}\in \mathcal H,$
\[
 \ \mathcal H \ne P^2(S | G_1,G_2,...,G_{N-1}).
 \]
In this case, $m \le N$ where $m$ is as in (1-1). 

Let $\mu$ be a compactly supported finite positive measure on the complex plane
$\mathbb {C}$ and let $spt(\mu)$ denote the support of $\mu .$ For a compact subset $K$ with $spt(\mu) \subset K,$ let $R^2(K, \mu)$ be the closure of $Rat(K)$ in $L^2(\mu).$
 Let $P^2(\mu)$ denote the closure of
$\mathcal{P}$ in $L^2(\mu).$ 

If $S$ is rationally cyclic, then $S$ is unitarily equivalent to multiplication by $z$ on $R^2(\sigma (S), \mu _1),$ where $m=1$ and $F_1 = 1.$ We may write $R^2(S | F_1) = R^2(\sigma (S), \mu _1).$ If $S$ is cyclic, then $S$ is unitarily equivalent to multiplication by $z$ on $P^2(\mu _1).$  We may write $P^2(S | F_1) = P^2(\mu _1).$

For a rationally $N-$cyclic subnormal operator $S$ with cyclic vectors $F_1,F_2,...,F_N$ and $\lambda \in \sigma (S),$  we denote the map
 \[
 \ E(\lambda ) : \sum_{i=1}^N r_i(S)F_i \rightarrow \left[ \begin{array}{c}
r_1(\lambda )\\
r_2(\lambda ) \\
... \\
r_N(\lambda ) \end{array} \right],\tag{1-7}
 \]
where $r_1,r_2,...,r_N\in Rat(\sigma(S)).$ If $E(\lambda)$ is bounded from $\mathcal K$ to $(\mathbb C^N, \|.\|_{1,N}),$ where $\|x\|_{1,N} = \sum_{i=1}^N |x_i|$ for $x\in \mathbb C^N,$
 then every component in the right hand side extends to a bounded linear functional on $\mathcal H$ and we will call $\lambda$ a
bounded point evaluation for $S.$ We use $bpe(S)$ to denote the set of bounded point evaluations for $S.$  The set $bpe(S)$ does not depend on the choices of cyclic vectors $F_1,F_2,...,F_N$ (see Corollary 1.1 in \cite{moz}).
A point $\lambda_0 \in int(bpe(S))$ is called an analytic bounded point evaluation for $S$ if there is a neighborhood $B(\lambda_0, \delta)\subset bpe(S)$ of $\lambda_0$ such that $E(\lambda)$ is analytic as a function of $\lambda$ on $B(\lambda_0, \delta)$ (equivalently (1-7) is uniformly bounded for  $\lambda\in B(\lambda_0, \delta)$).  We use $abpe(S)$ to denote the set of analytic bounded point evaluations for $S.$  The set $abpe(S)$ does not depend on the choices of cyclic vectors $F_1,F_2,...,F_N$ (also see Remark 3.1 in \cite{moz}). Similarly, for an $N-$cyclic subnormal operator $S,$ we can define $bpe(S)$ and $abpe(S)$ if we replace $r_1,r_2,...,r_N\in Rat(\sigma(S))$ in (1-7) by $p_1,p_2,...,p_N\in \mathcal P.$

\cite{bfp85} show that the Bergman shift has invariant subspaces with the codimension $N$ property for every $N\in \{1, 2,...,\infty \}.$ This means that on the Bergman space, the set of all analytic functions $f$ on the unit disk $\mathbb D$ satisfying
 \[
 \ \int_{\mathbb D} |f(z)|^2 d A(z) < \infty ,
 \]
where $A$ is area measure, there is a closed subspace $\mathcal M$ that is invariant under
multiplication by the independent variable $z,$ and such that
 \[
 \ dim \left ( \mathcal M / (z\mathcal M) \right ) = N.  
 \]
Their construction is abstract, and these subspaces are hard to envision.
Later, \cite{h93} gave a concrete construction, using zero sets
whose union is not a zero set. Moreover, \cite{ARS96} show that there are $f_1, f_2, ..., f_N$ such that 
 \[
 \ \mathcal M = P^2(S | f_1,...,f_N),
 \]
where $S$ is multiplication by $z$ operator on $\mathcal M,$ and $bpe(S) = abpe(S) = \mathbb D.$

For $N = 1,$ \cite{thomson} proves a remarkable structural theorem for $P^2(\mu ).$ 

\begin{TTheorem}
There is a Borel partition $\{\Delta_i\}_{i=0}^\infty$ of $spt\mu$ such that the space $P^2(\mu |_{\Delta_i})$ contains no nontrivial characteristic functions and
 \[
 \ P^2(\mu ) = L^2(\mu |_{\Delta_0})\oplus \left \{ \oplus _{i = 1}^\infty P^2(\mu |_{\Delta_i}) \right \}.
 \]
Furthermore, if $U_i$ is the open set of analytic bounded point evaluations for
$P^2(\mu |_{\Delta_i})$ for $i \ge 1,$ then $U_i$ is a simply connected region and the closure of $U_i$ contains $\Delta_i.$
\end{TTheorem}

\cite{ce93} extends some results of Thomson's Theorem to the space $R^2(K,\mu ),$ while \cite{b08} expresses $R^2(K,\mu )$ as a direct sum that includes both Thomson's theorem and results of \cite{ce93}. For a compactly supported complex Borel measure $\nu$ of $\mathbb C,$ by estimating analytic capacity of the set $\{\lambda: | C\nu (\lambda)| \ge c \},$ where $C\nu$ is the Cauchy transform of $\nu$ (see Section 2 for definition), \cite{b06}, \cite{ars}, and \cite{ARS10} provide interesting alternative proofs of Thomson's theorem. Both their proofs rely on X. Tolsa's deep results on analytic capacity. There are other related research papers for $N=1$ in the history. For example, \cite{b79}, \cite{h79}, \cite{bm}, and \cite{yang3}, etc.

Thomson's Theorem shows in 
Theorem 4.11 of \cite{thomson} that $abpe(S) = bpe(S)$ for a cyclic subnormal operator $S$ (See also Chap VIII Theorem 4.4 in \cite{conway}).
The results lead to the next question stated by Conway 7.11 p. 65 of \cite{conway}.

{\it Does $abpe(S) = Int(bpe(S))$ hold for an arbitrary rationally cyclic subnormal operator $S$?}

Corollary 5.2 in \cite{ce93} affirmatively answers the question. Our following theorem extends the result to rationally $N-$cyclic subnormal operators. 

\begin{Theorem}\label{MTheorem1} Let $S$ on $\mathcal H$ be a pure subnormal operator and let $M_z$ on $\mathcal K$ (satisfying (1-1) and (1-4)) be its minimal normal extension.

(1) If $S$ is $N-$cyclic, then $abpe(S)=bpe(S).$

(2) If $S$ is rationally $N-$cyclic, then $abpe(S)=bpe(S)\cap Int(\sigma(S)).$
\end{Theorem}

$T\in\mathcal L(\mathcal H )$ satisfies Bishop’s property ($\beta$) provided
that for any open set $O,$ and any sequence of analytic functions
$f_n : O \rightarrow \mathcal H,$ the convergence of $(T - \lambda)f_n(\lambda )$ to zero, on compact sets,
forces $f_n$ to converge to zero on compact sets. The class of operators
with Bishop’s property ($\beta$) is very large, in particular, it contains subnormal
operators (see \cite{ln2000}). There are some results related to Conway's question for rationally multicyclic operators on Banach spaces satisfying
Bishop’s property $(\beta).$ For example, \cite{moz} provides an example of a rationally multicyclic operator $T$ satisfying
Bishop’s property $(\beta),$ but $abpe(T) \ne Int(bpe(T)).$ Also \cite{mmn} studies rationally cyclic operators.

Applying Theorem \ref{MTheorem1}, we obtain the following results.

\begin{Corollary}\label{MCorollary}
Let $S$ on $\mathcal H$ be a pure subnormal operator and let $M_z$ on $\mathcal K$ (satisfying (1-1) and (1-4)) be its minimal normal extension. 

(1) Suppose $S$ is $N-$cyclic and $\lambda _0 \in \mathbb C.$ If $dim(Ker(S-\lambda _0 I)^* ) = N,$ then $Ran(S-\lambda _0 I)$ is closed and $\lambda_0 \in \sigma(S) \setminus \sigma_e (S).$

(2) Suppose $S$ is rationally $N-$cyclic and $\lambda _0 \in Int(\sigma (S)).$ If $dim(Ker(S-\lambda _0 I)^* ) = N,$ then $Ran(S-\lambda _0 I)$ is closed and $\lambda_0 \in \sigma(S) \setminus \sigma_e (S).$
\end{Corollary}

We prove Theorem \ref{MTheorem1} and Corollary \ref{MCorollary} in section 2.

\section{The Proofs}

Let $\nu$ be a compactly supported finite measure on $\mathbb {C}.$ The Cauchy transform
of $\nu$ is defined by
\[
\ C\nu (z) = \int \dfrac{1}{w - z} d\nu (w)
\]
for all $z\in\mathbb {C}$ for which
$\int \frac{d|\nu|(w)}{|w-z|} < \infty .$ A standard application of Fubini's
Theorem shows that $C\nu \in L^s_{loc}(\mathbb {C} )$ for $ 0 < s < 2,$ in particular, it is
defined for Area almost all $z,$ and clearly $C\nu$ is analytic in $\mathbb C_\infty \setminus spt \nu,$ where
$\mathbb C_\infty = \mathbb C \cup \{\infty \}$ is the Riemann sphere.

Now suppose that $\nu$ is a compactly supported finite measure on $\mathbb C$ that annihilates
the rational functions $Rat(spt(\nu )).$ Then, for $r\in Rat(spt(\nu )),$ 
\[
\ \int \dfrac{r(z) - r(w)}{z - w}d\nu(z ) = 0.
\]
Rearranging, we see that
\[
\ r(w) C\nu (w) = \int \dfrac{r(z)}{z - w} d\nu(z)
\]
for Area almost all $w.$ 

Suppose that $S$ on $\mathcal H$ is a pure rationally $N-$cyclic subnormal operator with cyclic vectors $F_1,F_2,...,F_N$ and let $M_z$ on $\mathcal K$ (satisfying (1-1) to (1-4)) be its minimal normal extension.
Let $G_i \in \mathcal K$ and $G_i \bot \mathcal H$ for $i = 1,2,...,N.$ Denote 
 \[
 \ \nu_{ij} = \left\langle F_i(z), G_j(z) \right\rangle \mu _1,\tag{2-1}
 \]
then $\nu _{ij}$ annihilates $Rat(\sigma(S)).$ We have the following estimation for $r_i\in Rat(\sigma(S))$ and $B(\lambda_0, \delta) \subset Int(\sigma(S)).$
\[
\begin{aligned}
\ \int _{B(\lambda_0, \delta)} \left | \sum_{i = 1}^N r_iC\nu_{ij}\right |dA & = \int _{B(\lambda_0, \delta)} \left | C ( \sum _{i = 1}^N r_i\nu_{ij} ) \right |dA \\
\ & \le \int \int _{B(\lambda_0, \delta)} \dfrac{1}{|z-w|} dA(z) \left | \sum_{i = 1}^N r_iF_i \right | |G_j| d\mu _1\\
\ & \le M \delta \left \| \sum_{i = 1}^N r_iF_i \right \|_{\mathcal K}, 
\end{aligned} \tag{2-2}
\]
where $M$ is a constant. Notice that $C$ refers to the Cauchy transform (not a constant).

For a compact $K \subset \mathbb C$ we
define the analytic capacity of $K$ by
\[
\ \gamma(K) = sup |f'(\infty)|
\]
where the sup is taken over those functions $f$ analytic in $\mathbb C_\infty \setminus K$ for which
$|f(z)| \le 1$ for all $z \in \mathbb C_\infty \setminus K,$ and
\[
\ f'(\infty) = \lim _{z \rightarrow \infty} z[f(z) - f(\infty)].
\]
The analytic capacity of a general $E \subset \mathbb C$ is defined to be 
\[
\ \gamma (E) = \sup \{\gamma (K) : K \subset E, ~K~ compact\}.
\]
Good sources for basic information about analytic
capacity are \cite{Gar72}, Chapter VIII of \cite{gamelin}, Chapter V of \cite{conway}, and \cite{tol14}.

A related capacity, $\gamma _+,$ is defined for $E \subset \mathbb C$ by
\[
\ \gamma_+(E) = sup \|\mu \|
\]
where now the sup is taken over positive measures $\mu$ with compact support
contained in $E$ for which $\|C\mu \|_{L^\infty (\mathbb C)} \le 1.$ Since $C\mu$ is analytic in $\mathbb C_\infty \setminus spt \mu$ and $(C \mu)'(\infty) = −\|\mu \|,$ we have
\[
\ \gamma _+(E) \le \gamma (E)
\]
for all $E \subset \mathbb C.$ \cite{Tol03} proves the astounding result (Tolsa's Theorem) that 
$\gamma_+$ and $\gamma$ are actually equivalent. 
 That is, there is an absolute constant $A_T$ such that
\[ 
\ \gamma (E) \le A_ T \gamma_+(E)
\]
for all $E \subset \mathbb C.$ The following semiadditivity of analytic capacity is a conclusion of Tolsa's Theorem.
\[
\ \gamma \left (\bigcup_{i = 1}^m E_i \right ) \le A_T \sum_{i=1}^m \gamma(E_i)\tag{2-3}
\]
where $E_1,E_2,...,E_m \subset \mathbb C.$

We set
 \[
 \ \sigma_0(S) = \begin{cases}\mathbb C, ~ \text{if $S$ is $N-$cyclic} \\ \sigma (S), \text{if $S$ is rationally $N-$cyclic} \end{cases} \tag{2-4}
\]

\begin{Theorem}\label{bpeTheorem} 
Let $S$ on $\mathcal H$ be a pure $N-$cyclic or rationally $N-$cyclic subnormal operator with cyclic vectors $F_1,F_2,...,F_N$ and let $M_z$ on $\mathcal K$ (satisfying (1-1) to (1-4)) be its minimal normal extension. Let $G_j \bot \mathcal H$ for $j=1,2,...,N$ and let $\nu_{ij}$ be as in (2-1). If $\lambda _0 \in Int(\sigma_0(S))$ and $\nu_{ij}$ satisfy

(1) 
\[
\ \int \dfrac{1}{|z - \lambda_0|} d |\nu_{ij} | ( z) < \infty 
\]
for all $i,j = 1,2,...,N.$ 

(2) The matrix
\[
\ \left [ C\nu_{ij} (\lambda _0) \right ] _{N\times N}
\]
is invertible.

Then $\lambda _0$ is an analytic bounded point evaluation for $S.$
\end{Theorem}

Before proving Theorem \ref{bpeTheorem}, let us use Theorem \ref{bpeTheorem} to prove Theorem \ref{MTheorem1} and Corollary \ref{MCorollary}. For a subnormal operator $S$ on $\mathcal H$ and its minimal normal extension $M_z$ on $\mathcal K$ (satisfying (1-1) to (1-4)) with $\mu_1(\{\lambda \}) > 0,$ we define
 \[
 \ \mathcal K _{\lambda} = \mathcal K |_{\{\lambda\}^c} =\{\chi_{\{\lambda\}^c}F: F\in \mathcal K\}, ~ M_z^{\lambda} = M_z |_{\mathcal K _{\lambda}}, 
 \] 
and
 \[
 \ \mathcal H _{\lambda} = \mathcal H |_{\{\lambda\}^c} = \{\chi_{\{\lambda\}^c}F: F\in \mathcal H\},~ S^{\lambda} = S |_{\mathcal H _{\lambda}}. 
 \]

\begin{Lemma}\label{removeAtom}
Let $S$ on $\mathcal H$ be a pure $N-$cyclic or rationally $N-$cyclic subnormal operator with cyclic vectors $F_1,F_2,...,F_N$ and let $M_z$ on $\mathcal K$ (satisfying (1-1) to (1-4)) be its minimal normal extension. Suppose $\mu_1(\{\lambda _0\}) > 0$ and $\lambda_0\in bpe(S).$
Define 
 \[
 \ T:~ F\in \mathcal H \rightarrow F^{\lambda_0} = \chi_{\{\lambda_0\}^c}F \in \mathcal H_{\lambda_0},
 \]
then $T$ is invertible and $S^{\lambda _0} = TST^{-1},$ that is, $S^{\lambda _0}$ is similar to $S.$
Consequently, $S^{\lambda _0}$ on $H _{\lambda_0}$ is a pure $N-$cyclic or rationally $N-$cyclic subnormal operator with cyclic vectors $F_1^{\lambda_0},F_2^{\lambda_0},...,F_N^{\lambda_0}$ and $\lambda_0\in bpe(S^{\lambda _0}).$
\end{Lemma}

Proof: Assume $S$ is pure and rationally $N-$cyclic (same proof for $N-$cyclic), $\mu _1(\{\lambda_0\}) > 0,$  and $\lambda_0\in bpe(S).$ Then there is a constant $M_0 >0$ such that
 \[
 \ \begin{aligned}
 \ |r_k(\lambda_0)|^2& \le  M_0  \left\| \sum_{i = 1}^N r_iF_i \right\|^2_{\mathcal K} \\
 \ & =  M_0 \left \| \sum_{i = 1}^N r_iF_i^{\lambda_0} \right \|^2_{\mathcal K_{\lambda_0}} + M_0 \left | \sum_{i = 1}^N r_i(\lambda_0) F_i (\lambda_0) \right |^2 \mu _1 (\{\lambda_0\}),
 \ \end{aligned} \tag{2-5}
\]
where $r_k\in Rat(\sigma(S))$ for $k = 1,2,...,N.$
Suppose that $\lambda_0$ is not a bounded point evaluation for $S^{\lambda _0},$ then there exist $N$ sequences of rational functions $\{r_{in}\}_{1\le i \le N, 1 \le n < \infty} \subset Rat(\sigma (S))$ such that 
 \[
 \ \left \| \sum_{i = 1}^N r_{in}F_i^{\lambda_0} \right \|^2_{\mathcal K_{\lambda_0}} \rightarrow 0
 \]
and $|r_{i_0 n} (\lambda_0 )|\rightarrow 1$ for some fixed $i_0.$ Set $a_n = |\sum_{i = 1}^N r_{in}(\lambda_0) F_i (\lambda_0)|,$ then from (2-5) for $k = i_0,$ there exists $c_0 > 0$ such that
 \[
 \ \liminf_{n\rightarrow\infty} a_n \ge c_0.
 \]
Let 
 \[
 \ H_n(z) = \dfrac{\sum_{i = 1}^N r_{in}(z) F_i (z)}{a_n},
 \]
then $H_n \in \mathcal H.$ By choosing a subsequence, we may assume there is $v\in \mathbb C^m$ and $v \ne 0$ such that  
 \[
 \ \|H_n - \chi_{\{\lambda_0\}}v \| _{\mathcal K}\rightarrow 0.
 \]
Therefore, $\chi_{\{\lambda_0\}} v \in \mathcal H$ and this is a contradiction since $S$ is pure. Hence, $\lambda_0\in bpe(S^{\lambda _0}).$
So there is a constant $M_1 >0$ such that
 \[
 \ |r_k(\lambda_0)|^2 \le M_1 \left \| \sum_{i = 1}^N r_iF_i^{\lambda_0} \right \|^2_{\mathcal K_{\lambda_0}},
\]
where $r_k\in Rat(\sigma(S))$ for $k = 1,2,...,N.$ Hence,
 \[
 \  \left\| \sum_{i = 1}^N r_iF_i \right\|^2_{\mathcal K} =  \left \| \sum_{i = 1}^N r_iF_i^{\lambda_0} \right \|^2_{\mathcal K_{\lambda_0}} +  \left | \sum_{i = 1}^N r_i(\lambda_0) F_i (\lambda_0) \right |^2 \mu _1 (\{\lambda_0\}) \le \left (1 + M_1( \sum_{i = 1}^N |F_i (\lambda_0) |)^2 \mu _1 (\{\lambda_0\})\right) \left \| \sum_{i = 1}^N r_iF_i ^{\lambda_0}\right \|^2_{\mathcal K_{\lambda_0}}.
\]
This implies that $T$ is invertible.

Proof of (2) in Theorem \ref{MTheorem1}: Suppose $\lambda_0\in bpe(S)\cap Int(\sigma (S)).$ By Lemma \ref{removeAtom}, we assume  that $\mu _1 (\{\lambda _0 \}) = 0.$ There are $g_1,g_2,...,g_N \in \mathcal H$ such that
 \[
 \ r_j(\lambda _0) = (\sum_{i = 1}^N r_i F_i, g_j),
 \]
where $r_k\in Rat(\sigma (S)).$ Set $G_j = \overline{z - \lambda _0} g_j$ for $j = 1,...,N.$ Let $\nu_{ij}$ be as in (2-1), then
 \[
\ \int \dfrac{1}{|z - \lambda_0|} d |\nu_{ij} | ( z) \le \|F_i\|\|g_j\|< \infty , 
\]
and
 \[
 \ C\nu_{ij} (\lambda _0) = (F_i,g_j) = \delta_{ij}.
 \]
By Theorem \ref{bpeTheorem}, we conclude $\lambda_0\in abpe(S).$

The proof of (1) in Theorem \ref{MTheorem1} is the same.

Proof of (2) in Corollary \ref{MCorollary}: From the assumptions of the corollary  and Theorem \ref{MTheorem1}, we see $\lambda_0\in abpe(S).$ There are $\delta , M > 0$ such that 
 \[
 \ |r_j(\lambda)| \le  M \left\| \sum_{i = 1}^N r_iF_i \right\|^2_{\mathcal K}
 \] 
for $B(\lambda_0,\delta)\subset Int(\sigma (S)),$ $\lambda\in B(\lambda_0,\delta),$ and $r_j\in Rat(\sigma (S)).$
Using the maximal modulus principle,
\[
 \ \sup _{1 \le j\le N,\lambda\in B(\lambda_0,\delta)} \left | r_j(\lambda ) \right | \le\dfrac{M}{\delta} \left \|(z-\lambda_0) \sum_{i = 1}^N r_iF_i \right \|_{\mathcal K}.
\]
Hence,
\[
 \ \int \left |\sum_{i = 1}^N r_iF_i\right |^2 d\mu _1 \le \int_{B(\lambda_0,\delta)^c} \left |\sum_{i = 1}^N r_iF_i \right |^2 d\mu _1 + (\sum_{i = 1}^N \|F_i\| ) ^2\sup _{1 \le j\le N,\lambda\in B(\lambda_0,\delta)} \left | r_j(\lambda ) \right |^2.
\]
Therefore,
\[
 \ \left\| \sum_{i = 1}^N r_iF_i \right \|_{\mathcal K} \le M_1 \left\|(z - \lambda _0) \sum_{i = 1}^N r_iF_i \right\|_{\mathcal K},
\]
where
 \[
 \ M_1^2 = \left ( 1  + M^2 (\sum_{j =1}^N \| F_j\|)^2 \right )/ \delta^2.
\]
So $Ran(S - \lambda _0)$ is closed. The corollary is proved.

The proof of (1) in Corollary \ref{MCorollary} is the same.

Theorem \ref{bpeTheorem} is a generalization of Corollary 2.2 in \cite{ars} where $N = 1.$ There are fundamental differences between $N=1$, where the existence of analytic bounded point evaluations for $P^t(\mu )$ was first proved in \cite{thomson}, and $N > 1,$ where analytic bounded point evaluations may not exist (see  the example at the end of this section). To prove Theorem \ref{bpeTheorem}, we need several lemmas.

The following Lemma is from Lemma B in \cite{ars}.

\begin{Lemma} \label{lemmaARS}
There are absolute constants $\epsilon _1 > 0$ and $C_1 < \infty$ with the
following property. For $R > 0,$ let $E \subset clos (R\mathbb D)$ with 
$\gamma (E) < R\epsilon_1.$ Then
\[
\ |p(0)| \le \dfrac{C_1}{R^2} \int _{clos ( R\mathbb D)\setminus E} |p| \frac{dA}{\pi}
\]
for all $p \in \mathcal P.$
\end{Lemma}

\begin{Lemma} \label{lemmaARSEx}
Let $\epsilon _1 > 0$ and $C_1 < \infty$ be as in Lemma \ref{lemmaARS}. For $R > 0,$ let $E \subset clos (R\mathbb D)$ with 
$\gamma (E) < \frac{R}{2}\epsilon_1.$ Then
\[
\ |p(\lambda)| \le \dfrac{4C_1}{R^2} \int _{clos ( R\mathbb D)\setminus E} |p| \frac{dA}{\pi}
\]
for all $\lambda\in B(0, \frac{R}{2})$ and $p \in \mathcal P.$
\end{Lemma}

Proof: For $\lambda\in B(0, \frac{R}{2}),$ let $E_{\lambda} =  \bar B(\lambda, \frac{R}{2})\cap E - \lambda .$ Then $E_{\lambda} \subset clos (\frac{R}{2}\mathbb D)$ and $\gamma (E_{\lambda}) \le \gamma (E) < \frac{R}{2}\epsilon_1.$ From Lemma \ref{lemmaARS}, we have
\[
\ |p(0)| \le \dfrac{C_1}{(\frac{R}{2})^2} \int _{clos ( \frac{R}{2}\mathbb D)\setminus E_{\lambda}} |p| \frac{dA}{\pi}.
\]
Replacing $p(z)$ by $p(z+\lambda), $ we get 
\[
\ |p(\lambda)| \le \dfrac{4C_1}{R^2} \int _{clos ( \frac{R}{2}\mathbb D)\setminus E_{\lambda}} |p(z+\lambda)| \frac{dA(z)}{\pi} = \dfrac{4C_1}{R^2} \int _{\bar B(\lambda, \frac{R}{2})\setminus E} |p| \frac{dA}{\pi}\le \dfrac{4C_1}{R^2} \int _{clos ( R\mathbb D)\setminus E} |p| \frac{dA}{\pi}.
\]

Let $\nu$ be a compactly supported finite measure on $\mathbb {C}.$ For $\epsilon > 0,$ $C_\epsilon \nu$ is defined by
\[
\ C_\epsilon \nu (z) = \int _{|w-z| > \epsilon}\dfrac{1}{w - z} d\nu (w),
\]
and the maximal Cauchy transform is defined by
 \[
 \ C_* \nu (z) = \sup _{\epsilon > 0}| C_\epsilon \nu (z) |.
 \]
From Proposition 2.1 of \cite{Tol02} and Tolsa's Theorem, we have the following estimation (also see \cite{tol14} Proposition 4.16):
\[
\ \gamma(\{|C_*\nu | \geq a\}) \le \dfrac{C_T}{a} \|\nu \|, \tag{2-5}
\]
 where $C_T$ is an absolute positive constant and $a > 0.$

\begin{Lemma}\label{CauchyTLemma} 
Suppose $\nu$ is a finite compactly supported Borel measure of $\mathbb C,$ $\lambda _0\in \mathbb C,$ and 
 \[
 \   \int \dfrac{1}{|z - \lambda _0| } d| \nu | (z) < \infty . \tag{2-6}
 \] 
If $\epsilon_0 , a > 0,$ then 

(1)  there exist $0 < \delta_a < \frac{1}{4}$ such that 
\[
 \ 2\sqrt{\delta_a} \int \dfrac{1}{|z - \lambda _0| } d| \nu | (z) + \int _{B (\lambda _0, 2 \delta _a)} \dfrac{1}{|z - \lambda _0| } d| \nu | (z) < \frac{a}{2}\tag{2-7}
 \]
and
 \[
 \ \dfrac{2C_T}{a} \int _{B (\lambda _0, \sqrt{ \delta _a})} \dfrac{1}{|z - \lambda _0| } d| \nu | < \epsilon_0; \tag{2-8}
 \]

(2) for  $0 < \delta < \delta_a,$ there exists $E_\delta \subset \bar B (\lambda _0, \delta )$ such that $\gamma(E_\delta) <\epsilon _0 \delta$ and 
\[
\ |C\nu (\lambda) - C\nu (\lambda _0)  | \le a
\]
almost everywhere with respect to the area measure on $B (\lambda _0, \delta ) \setminus E_\delta .$
\end{Lemma}

Proof: (2-7) and (2-8) of (1) follow from (2-6). For (2), we fix $0 < \delta < \delta_a.$ 
Let $\nu_\delta = \frac{\chi _{B (\lambda _0, \sqrt{ \delta })}}{z - \lambda _0 } \nu .$  For $\epsilon < \delta $ and $\lambda \in B (\lambda _0, \delta ),$ we get:
 \[
 \ \bar B (\lambda , \epsilon) \subset B (\lambda _0, 2 \delta ) \subset B (\lambda _0, \sqrt{\delta}),
 \]
 \[
 \ (\bar B (\lambda , \epsilon))^c \cap B (\lambda _0, \sqrt{\delta} )^c = B (\lambda _0, \sqrt{\delta} )^c \subset B (\lambda , \sqrt{\delta} - \delta )^c,
 \]
and
 \[ 
 \ \begin{aligned}
 \ & |C_\epsilon \nu (\lambda) - C\nu (\lambda _0)| \\
 \ \le & |\lambda - \lambda _0| \left | \int _{|z - \lambda| > \epsilon}\dfrac{d\nu}{(z - \lambda)(z - \lambda_0)} \right | + \int _{\bar B (\lambda, \epsilon)} \dfrac{1}{|z - \lambda _0| } d| \nu | (z) \\
\ \le & \delta  \left | \int _{(\bar B (\lambda , \epsilon))^c \cap B (\lambda _0, \sqrt{\delta} )^c}\dfrac{d\nu }{(z - \lambda)(z - \lambda_0)} \right | + \delta \left |\int _{|z - \lambda| > \epsilon}\dfrac{d\nu_\delta}{(z - \lambda)} \right | + \int _{B (\lambda _0, 2 \delta  )} \dfrac{1}{|z - \lambda _0| } d| \nu | (z) \\
 \ \le & \delta\int _{|z - \lambda | \ge \sqrt{\delta} - \delta } \dfrac{1}{|z - \lambda ||z - \lambda _0| } d| \nu | (z) + \delta |C_\epsilon \nu _\delta (\lambda )|  + \int _{B (\lambda _0, 2 \delta )} \dfrac{1}{|z - \lambda _0| } d| \nu | (z) \\
 \ \le & 2\sqrt{\delta}\int \dfrac{1}{|z - \lambda _0| } d| \nu | (z) + \delta C_* \nu _\delta (\lambda )  + \int _{B (\lambda _0, 2 \delta )} \dfrac{1}{|z - \lambda _0| } d| \nu | (z) \\
 \ \le & \dfrac{a}{2} + \delta  C_* \nu _\delta (\lambda ), 
 \ \end{aligned} \tag{2-9}
 \]
where the last two steps follow from $\frac{\delta }{\sqrt{\delta} - \delta } \le 2 \sqrt{\delta}$ and (2-7).
Let
 \[
 \ E_\delta= \{\lambda : C_* \nu _\delta (\lambda ) \ge \frac{a}{2\delta} \} \cap \bar B (\lambda _0, \delta ),
 \]
then from (2-9), we get
 \[
 \ \{\lambda : |C_\epsilon \nu (\lambda) - C\nu (\lambda _0)| \ge a \} \cap \bar B (\lambda _0, \delta ) \subset E_\delta.
 \]
From (2-5) and (2-8), we get
 \[
 \ \gamma (E_\delta) \le \dfrac{2C_T\delta}{a} \| \nu _\delta \| < \epsilon_0 \delta.
 \]
On $B (\lambda _0, \delta ) \setminus E_\delta,$ for $\epsilon < \delta ,$ we conclude that
 \[
 \ |C_\epsilon \nu (\lambda) - C\nu (\lambda _0)| < a.
 \]
The lemma follows since
 \[
 \ \lim_{_\epsilon\rightarrow 0} C_\epsilon \nu (\lambda) = C\nu (\lambda )~ a.e. ~ Area.
 \]

\begin{Remark}
(1) In \cite{ars} and \cite{ARS10}, a key step for their alternative proofs of Thomson's theorem is to show that $|C\nu (\lambda) |$ is bounded below on $B (\lambda _0, \delta ) \setminus E,$ where $\gamma (E) < \epsilon _0 \delta$ and $C\nu (\lambda _0) \ne 0.$ This is directly implied by above lemma. So the lemma provides an alternative proof of the property.

(2) For $\nu ,$ $\lambda_0 ,$ and $a > 0$ in Lemma \ref{CauchyTLemma}, we have
 \[
 \ \begin{aligned}
 \ & Area (\{ |C\nu (\lambda) - C\nu (\lambda _0)  | > a \} \cap B(\lambda _0, \frac{1}{n})) \\
 \ \le & \dfrac{1}{a} \int _{B(\lambda _0, \frac{1}{n})} |C\nu (\lambda) - C\nu (\lambda _0)  | dA (\lambda ) \\ 
 \ \le & \dfrac{1}{a} \int _{B(\lambda _0, \frac{1}{n})} |\lambda - \lambda _0| \left | C (\dfrac{\nu}{z-\lambda _0})(\lambda )\right | dA(\lambda ).
 \ \end{aligned}
 \] 
Therefore, by Lemma 1 of \cite{b67}, we see that $| C\nu (\lambda) - C\nu (\lambda_0)| \le a$ on a set having full area density at $\lambda_0$ whenever $|\lambda - \lambda_0 |$
is sufficiently small. Lemma \ref{CauchyTLemma} (2) shows that this inequality holds capacitary density which is needed (not area density as Browder considers) in order to apply Lemma \ref{lemmaARSEx} in proving our main theorem below.
\end{Remark}

The following lemma is a simple linear algebra exercise.

\begin{Lemma}\label{basicLA} 
Let $\|.\|_{1,n}$ be the $l^1$ norm of $\mathbb C^n,$ that is, $\|x\|_{1,n} = \sum_{i = 1}^n |x_i|.$ Let $A=[a_{ij}]_{n\times n}$ be an $n\times n$ matrix such that $|a_{ij} - \delta_{ij}| \le \frac{1}{2n}.$ Then
\[
\ \|xA\|_{1,n} \geq \epsilon_2 \|x\|_{1,n}, ~ x\in \mathbb C^n,
\]
where $\epsilon_2 > 0$ is an absolute constant.
\end{Lemma}

Proof of Theorem \ref{bpeTheorem}: 
Assume $S$ is rationally $N-$cyclic and $\lambda_0 =0.$
Let $A = [a_{ij} ]_{N\times N}$ be the inverse of the matrix in condition (2). By replacing $F_i$ by $\sum_{k = 1}^N a_{ik}F_k,$ we may assume that the matrix in the condition (2) is identity.

For $a= \frac{1}{2N},$ let $\delta_a^{ij}$ be $\delta_a$ for $ \nu  = \nu _{ij}$ in (1) of Lemma \ref{CauchyTLemma}. Set
 \[
 \ \delta _0 = \min_{1\le i,j \le N} \delta_a^{ij}.
 \] 
From (2) of Lemma \ref{CauchyTLemma}, for a given $0 < \delta <\delta _0$ with $B (0, \delta )\subset Int(\sigma(S))$ and $\epsilon_0 = \frac{\epsilon_1}{2A_TN^2},$ where $\epsilon_1$ is in Lemma \ref{lemmaARS} and $A_T$ is from (2-3), there exists $E_\delta ^{ij}\subset \bar B (0, \delta )$ such that $\gamma(E_\delta ^{ij}) <\epsilon _0 \delta$ and 
\[
\ |C\nu _{ij}(\lambda) - \delta _{ij}  | = |C\nu _{ij}(\lambda) - C\nu _{ij}(0)  | \le a
\]
almost everywhere with respect to the area measure on $B (0, \delta ) \setminus E_\delta ^{ij}.$ Set $E = \cup_{i,j =1}^N E_\delta ^{ij},$ $B (\lambda) = [C\nu _{ij}(\lambda)]_{N\times N},$ $R(\lambda ) = (r_1(\lambda ), r_2(\lambda ),...,r_N(\lambda ))$ where $r_1, r_2,...,r_N \in Rat(\sigma (S)),$ and $R(\lambda )B(\lambda ) = (b_1(\lambda ), b_2(\lambda ),...,b_N(\lambda )),$ then, from Lemma \ref{basicLA}, we have
 \[
 \ \sum_{j = 1}^N |b_j(\lambda )| = \|R(\lambda )B(\lambda ) \|_{1,N} \ge \epsilon_2 \|R(\lambda ) \|_{1,N} = \epsilon_2 \sum_{j = 1}^N |r_j(\lambda )|,
 \]
almost everywhere with respect to the area measure on $B (0, \delta ) \setminus E,$
 where
 \[
 \ b_j(\lambda ) = \sum_{i = 1}^N r_i (\lambda ) C(\nu _{ij})(\lambda ).
 \]
From (2-2), we see that
 \[
 \ \int _{B(0,\delta)}\sum_{j = 1}^N |b_j(\lambda )| dA(\lambda )\le NM\delta \left \|\sum_{i = 1}^N r_iF_i \right \|.
 \]
On the other hand, from (2-3), we see 
 \[
 \ \gamma (E) \le A_T \sum_{i,j = 1}^N \gamma (E_\delta ^{ij}) < A_TN^2\epsilon_0\delta = \epsilon_1\frac{\delta}{2} ,
 \]
Therefore, applying Lemma \ref{lemmaARSEx} for $r_1, r_2,...,r_N \in Rat(\sigma (S))$ since $r_1, r_2,...,r_N$ are analytic on $B(0,\delta),$ we conclude
\[
\begin{aligned}
\ \sum_{i = 1}^N |r_i (\eta) | &\le \dfrac{4C_1}{\delta^2} \int _{B(0,\delta)\setminus E} \sum_{i = 1}^N |r_i| \frac{dA}{\pi} \\
\ & \le \dfrac{4C_1}{\pi\epsilon_2 \delta^2} \int _{B(0,\delta)\setminus E} \sum_{j = 1}^N |b_j(\lambda )|  dA(\lambda ) \\
\ & \le \dfrac{4C_1NM}{\pi\epsilon_2 \delta} \left  \| \sum_{i = 1}^N r_i F_i \right \|
\end{aligned}
\] 
for $\eta \in B(0, \frac{\delta}{2}).$ So $0$ is an analytic bounded point evaluation for $S$.

If $S$ is $N-$cyclic, we just need to drop the condition that $B(0, \delta) \subset Int(\sigma(S))$ and replace  rational functions $r_1,r_2,...,r_N\in Rat(\sigma(S))$ by polynomials $p_1,p_2,...,p_N\in \mathcal P$ in the proof.
This completes the proof.
 
\begin{Example} A Swiss cheese $K$ can be constructed as
 \[
 \ K = \bar {\mathbb D} \setminus \cup_{n=1}^\infty B(a_n, r_n),
 \] 
where $B(a_n, r_n) \subset \mathbb D,$ $\bar B(a_i, r_i) \cap \bar B(a_j, r_j) = \emptyset $ for $i\ne j,$ $\sum_{n=1}^\infty r_n< \infty ,$ and $K$ has no interior points. Let $\mu$ be the sum of the arc length measures of $\partial \mathbb D$ and  all $\partial B(a_n, r_n).$ Let $\nu$ be the sum of $dz$ on  $\partial \mathbb D$ and all $-dz$ on $\partial B(a_n, r_n).$ For a rational function $f$ with poles off $K,$ we have
 \[
 \ \int f d\nu = 0.
 \]
Clearly $| \frac{d\nu}{d\mu} | > 0, ~ a.e. ~\mu$ and $\overline{(\frac{d\nu}{d\mu})} \perp R^2(K, \mu),$ so $S_\mu$ (multiplication by $z$)  on $R^2(K, \mu)$ is a pure subnormal operator and $\sigma (S_\mu ) = \sigma _e (S_\mu ) = K.$ Moreover, there exists a function $F \in R^2(K, \mu)$ such that $R^2(K, \mu) = P^2(S_\mu| 1, F).$
\end{Example}

Proof: Let 
 \[
 \ H_n (z) = \left ( \dfrac{r_1}{z-a_1}\right )^n \left ( \dfrac{r_2}{z-a_2}\right )^{n-1}...\left ( \dfrac{r_{n-1}}{z-a_{n-1}}\right )^2 \left ( \dfrac{r_n}{z-a_n}\right ).
 \]
Let $\{b_n\}$ be a sequence of positive numbers satisfying 
 \[
 \ \lim_{n\rightarrow\infty} \sum_{k = n+1}^\infty \dfrac{b_k}{b_n} = 0.
 \]
 Set $F=\sum_{n=1}^\infty b_n H_n.$ By construction, we see that for $z\in K,$ $|F(z)| \le \sum_{n=1}^\infty b_n < \infty .$ Let
 \[
 \ p_1 = - \dfrac{r_1}{b_n(z-a_1)H_n} \sum_{k=1}^{n-1} b_k H_k, ~ p_2 =  \dfrac{r_1}{b_n(z-a_1)H_n}. 
 \]
Then $p_1$ and $p_2$ are polynomials and for $z\in K,$
 \[
 \ \left | p_1 + p_2 F - \dfrac{r_1}{z-a_1} \right | \le  \sum_{k = n+1}^\infty \dfrac{b_k}{b_n} \rightarrow 0.
 \] 
Hence $\frac{1}{z-a_1} \in P^2(S_\mu| 1, F).$ Similarly, one can prove that $\frac{1}{(z-a_n)^m} \in P^2(S_\mu| 1, F).$ Therefore,  $R^2(K, \mu) = P^2(S_\mu| 1, F),$ rationally cyclic subnormal operator $S_\mu$ is $2-$cyclic, and $abpe(S) = \emptyset.$
        
\vspace{\baselineskip}
\vspace{\baselineskip}
\vspace{\baselineskip}
      
\centerline{\bf Acknowledgment}

The author would like to thank the referee for providing very helpful comments and pointing out the fact in Remark (2).

\bibliography{Bibliography}

\end{document}